\documentclass[12pt,twoside,reqno]{amsart}

\usepackage{amsmath,amssymb,amsthm, amsfonts }
\usepackage{fancyhdr}
\usepackage{epsfig,color,colordvi}

\usepackage[numbers,sort&compress,square,comma]{natbib}

\theoremstyle{definition}

\theoremstyle{remark}

\theoremstyle{plain}

\newcommand{\be}{\begin{equation}}
\newcommand{\ee}{\end{equation}}
\newcommand{\nn}{\nonumber}

\def\LL{\mathbf{L}}

\def\ee{\mathbb{E}}
\def\nn{\mathbb{N}}
\def\rr{\mathbb{R}}
\def\pp{\mathbb{P}}

\def\LL{\mathcal L}

\newtheorem{thm}{Theorem}[section]
\newtheorem{cor}[thm]{Corollary}
\newtheorem{lem}[thm]{Lemma}
\newtheorem{prop}[thm]{Proposition}
\newtheorem{remarks}[thm]{Remarks}

\def\bdes{\begin{description}}
\def\ndes{\end{description}}
\def\beq{\begin{equation}}
\def\deq{\end{equation}}

\def\bdef{\begin{defn}}
\def\ndef{\end{defn}}

\def\bthm{\begin{thm}}
\def\nthm{\end{thm}}

\def\bprop{\begin{prop}}
\def\nprop{\end{prop}}

\def\brmk{\begin{remarks}}
\def\nrmk{\end{remarks}}

\def\bexa{\begin{exa}}
\def\nexa{\end{exa}}

\def\blem{\begin{lem}}
\def\nlem{\end{lem}}

\def\bcor{\begin{cor}}
\def\ncor{\end{cor}}

\def\dsp{\displaystyle}

\def\bexe{\begin{exe}}
\def\nexe{\end{exe}}

\def\bprf{\begin{proof}}
\def\nprf{\end{proof}}

\def\fQ{{{\rm Q}\kern-.65em {}^{{}_/ }\,}}
\def\fQQ{ {{\rm Q}\kern-.57em \scriptscriptstyle{}^{]\kern.055em[}\,}}

\def\ord{\kern0.1em o\kern-0.02em{}_{\ds\breve{}}\kern0.1em}
\def\Ord{\kern0.1em O\kern-0.02em{\ds\breve{}}\kern0.1em}
\def\ds{\displaystyle}

\def\fmonth{\ifcase\month\or Jan\or Feb\or Mar\or Apr
\or May\or Jun\or Jul\or Aug\or Sep
\or Oct\or Nov\or Dec\fi\ }
\def\mmddyyyy{\the\month.\the\day.\the\year}
\def\ddmmyyyy{\the\day.\the\month.\the\year}
\def\Mddyyyy{\fmonth~\the\day,~\the\year}

\providecommand{\pp}[1]{\langle#1\rangle}

\numberwithin{equation}{section}

\allowdisplaybreaks[1]

\begin{document}


\author{Yu Miao}   
\address{Department of Mathematics and Statistics, Wuhan
University, 430072 Hubei, China and College of Mathematics and
Information Science, Henan Normal University, 453007 Henan, China.}
\email{yumiao728@yahoo.com.cn}
\author{Guangyu Yang}  
\address{Department of Mathematics and Statistics,
Wuhan University, 430072 Hubei, China.}
\email{study\_yang@yahoo.com.cn}

\date{May 26, 2006}






\begin{abstract}
 H\"ormann (2006) gave an extension of almost sure
central limit theorem for bounded Lipschitz $1$ function. In this
paper, we show that his result of almost sure central limit theorem
is also hold for any Lipschitz function under stronger conditions.

\end{abstract}




\title[A note for extension of almost sure central limit theory]
{A note for extension of almost sure central limit theory}

\maketitle


\def\X{\mathcal X}
\def\B{\mathcal B}

\section{Introduction}

The classical results on the almost sure central limit theorem
(ASCLT) dealt with partial sums of random variables. A general
pattern is that, if $X_1, X_2, \ldots$ be a sequence of independent
random variables with partial sums $S_n=X_1+\cdots+X_n$ satisfying
$(S_n-b_n)/a_n\xrightarrow{\LL}H$ for some sequences $a_n>0$,
$b_n\in\rr$ and some distribution function $H$, then under some mild
conditions we have
$$\lim_{n\to\infty}\frac{1}{\log n}\sum_{k=1}^n\frac{1}{k}I\Big\{(S_k-b_k)/a_k\leq x\Big\}=H(x)\ \ a.s.$$
for any continuity point $x$ of $H$.

Several papers have dealt with logarithmic limit theorems of this
kind and the above relation has been extended in various directions.
Fahrner and Stadtm\"uller \cite{FS} gave an almost sure version of a
maximum limit theorem. Berkes and Horv\'ath \cite{BH} obtained a
strong approximation for the logarithmic average of sample extremes.
 Berkes and Cs\'aki \cite{BC} showed that not only the central limit theorem, but every
weak limit theorem for independent random variables has an analogous
almost sure version. For stationary Gaussian sequences with
covariance $r_n$, Cs\'aki and Gonchigdanzan \cite{CG} proved an
almost sure limit theorem for the maxima of the sequences under the
condition $r_n\log n(\log\log n)^{1+\varepsilon}=O(1)$. For some
dependent random variables, Peligrad and Shao \cite{PS} and
Dudzi\'nski \cite{Du} obtained corresponding results about the
almost sure central limit theorem.

Recently, H\"ormann \cite{HS} gave an extension of almost sure
central limit theory under some regularity condition as the
following form:
 \beq\label{110}\lim_{N\to\infty}D^{-1}_N\sum_{k=1}^N
d_k f\Big(\frac{S_k}{a_k}-b_k\Big)=\int_{-\infty}^{\infty}f(x)dH(x)\
\ a.s.
 \deq
  where $f$ is a bounded Lipschitz $1$ function and
$D_N=\sum_{k=1}^Nd_k$, $\{d_k\}_{k\geq 1}$ is a sequence of positive
constants. Using his method, we will show that for any Lipschitz
function $f$, (\ref{110}) holds under some additional conditions.

At first, we give our main result.
\begin{thm}\label{thm1} Let $X_1, X_2, \ldots$ be independent random variables
with partial sums $S_n$ and assume that
 \bdes
 \item{$(C_1)$} For some numerical sequences $a_n>0$ and $b_n$, we have
   $\dsp\frac{S_n}{a_n}-b_n\xrightarrow{\LL}H$,
   where $H$ is some (possibly degenerate) distribution function.

 \item{$(C_2)$}
 $kd_k=O(1)$ and for some $0<\alpha<1$, $d_kk^\alpha$
 is eventually non-increasing.

 \item{$(C_3)$} For some $\rho>0$, $\dsp d_k=O\Big(\frac{D_k}{k(\log
 D_k)^\rho}\Big)$.

\item{$(C_4)$} There exist positive constants $C$, $\beta$, such
that
 $ a_k/a_l\leq
C(k/l)^\beta\ \ (1\leq k\leq l).$
 Furthermore, for some $0<r<\rho$,
  \beq\label{111}
\ee\Big|\frac{S_n}{a_n}-b_n\Big|^\mu=O(1),\;\ \ for\; some\;
integer\; \mu\geq(2\vee 4/(\rho-r)).
 \deq

 \ndes
 Then for any Lipschitz function $f$ on the real line, we have
 (\ref{110}).
\end{thm}
\begin{remarks}
{\rm Obviously, for any bounded Lipschitz 1 function $f$, under the
above assumptions, the equation (\ref{110}) holds, i.e. we can
obtain Theorem 1 in \cite{HS}.}

\end{remarks}
\section{Proof of Theorem \ref{thm1}}

In this section, we will give the proof of Theorem \ref{thm1},
according to the process of H\"ormann in \cite{HS}.

\begin{lem}\label{lem1} (See Lemma 1 in \cite{HS} )
Let $(D_{N})$ be a summation procedure, then the condition $(C_3)$
of Theorem \ref{thm1} implies that $D_N=o(N^\varepsilon)$ for any
$\varepsilon>0$.
\end{lem}

\begin{lem}\label{lem2}
Assume that condition $(C_4)$ of Theorem \ref{thm1} is satisfied and
$b_n=0$. Then for every Lipschitz function $f: \rr\to\rr$ there
exists constant $c>0$ such that
 \beq\label{210}
 \Big|Cov\Big(f\Big(\frac{S_k}{a_k}\Big),f\Big(\frac{S_l}{a_l}\Big)\Big)\Big|\leq
 c(k/l)^\beta
 \ \ (1\leq k\leq l),
 \deq
where $\beta$ is the same as in $(C_4)$.
\end{lem}

\begin{proof}
Firstly, we assume $f(0)=0$. Denoting $\|f\|$ the Lipschitz constant
of $f$, we get, by using the independence of $S_k$ and $S_l-S_k$,
 $$\aligned
 &\Big|Cov\Big(f\Big(\frac{S_k}{a_k}\Big),f\Big(\frac{S_l}{a_l}\Big)\Big)\Big|
 =
 \Big|Cov\Big(f\Big(\frac{S_k}{a_k}\Big),f\Big(\frac{S_l}{a_l}\Big)-f\Big(\frac{S_l-S_k}{a_l}\Big)\Big)\Big|\\
\leq&
\ee\Big|f\Big(\frac{S_k}{a_k}\Big)\Big[f\Big(\frac{S_l}{a_l}\Big)-f\Big(\frac{S_l-S_k}{a_l}\Big)\Big]\Big|
+\ee\Big|f\Big(\frac{S_k}{a_k}\Big)\Big|\ee\Big|f\Big(\frac{S_l}{a_l}\Big)-f\Big(\frac{S_l-S_k}{a_l}\Big)\Big|\\
\leq &
\|f\|^2\frac{a_k}{a_l}\ee\Big[\frac{S_k^2}{a_k^2}\Big]+\|f\|^2\frac{a_k}{a_l}(\ee\Big[\frac{S_k}{a_k}\Big])^2\\
\leq&2\|f\|^2\frac{a_k}{a_l}\ee\Big[\frac{S_k^2}{a_k^2}\Big]\leq
2C\|f\|^2\ee\Big[\frac{S_k^2}{a_k^2}\Big](k/l)^\beta,
 \endaligned$$
 where the last inequality is due to condition $(C_4)$. Since the equation (\ref{111}),
 we can take $c>0$ such that for any $k\geq 1$,
 $\dsp 2C\|f\|^2\ee\Big[\frac{S_k^2}{a_k^2}\Big]\leq c$. So (\ref{210}) is
 obtained.

If $f(0)\neq0$, we can define a function $g$, such that $
g(x)=f(x)-f(0)$, then $g$ is a Lipschitz function and $g(0)=0$. And
noting that,
\begin{align*}
Cov\Big(f\Big(\frac{S_k}{a_k}\Big),f\Big(\frac{S_l}{a_l}\Big)\Big)
 =Cov\Big(g\Big(\frac{S_k}{a_k}\Big),g\Big(\frac{S_l}{a_l}\Big)\Big)
\end{align*}
we complete the proof of the lemma.

\end{proof}
\begin{remarks}\label{rm1}
{\rm It is obvious to see that if we replace $\beta$ by any
$0<\beta^{'}<\beta$, the Lemma \ref{lem2} also holds. Hence, without
loss of generality, we can assume that $\beta$ is the same as
$\alpha$ in condition $(C_2)$ of Theorem \ref{thm1}.}
\end{remarks}
Next we will use the following notations,
 \beq\label{212}
 \xi_l:=f\Big(\frac{S_l}{a_l}\Big)-\ee f\Big(\frac{S_l}{a_l}\Big), \
 \ \xi_{k,l}:=f\Big(\frac{S_l-S_k}{a_l}\Big)-\ee
 f\Big(\frac{S_l-S_k}{a_l}\Big).
 \deq
\begin{lem}\label{lem3}
Assume that condition $(C_4)$ of Theorem \ref{thm1} is satisfied and
$b_n=0$, and define $\xi_l$ and $\xi_{k,l}$ as in (\ref{212}). If
$\{d_k, k\geq 1\}$ are arbitrary positive weights, then we have for
any $k\leq m\leq n$ and $p\in\nn$, $p\leq \mu$,
$$\ee\Big|\sum_{l=m}^nd_l(\xi_l-\xi_{k,l})\Big|^p\leq E_p\Big(\sum_{l=m}^nd_l^2l\Big)^{p/2},$$
where $\dsp
E_p=c_pC^{p}\|f\|^p\Big[\Big(\frac{2^{\kappa}}{\kappa}\Big)\vee\Big(1+\frac{1}{\kappa}\Big)\Big]^{p/2}$
and $\kappa=2\beta.$
\end{lem}
\begin{proof} Let $Q(l)=Q(k,l)=\xi_l-\xi_{k,l}$, then
$$\aligned
\ee|Q(l)|^p =&
\ee\Big|f\Big(\frac{S_l}{a_l}\Big)-f\Big(\frac{S_l-S_k}{a_l}\Big)-
\ee\Big[f\Big(\frac{S_l}{a_l}\Big)-f\Big(\frac{S_l-S_k}{a_l}\Big)\Big]\Big|^p\\
\leq
&\|f\|^p(a_k/a_l)^p\ee\Big(\frac{|S_k|}{a_k}+\ee\Big(\frac{|S_k|}{a_k}\Big)\Big)^p\\
\leq&
C^{p}\|f\|^p\ee\Big(\frac{|S_k|}{a_k}+\ee\Big(\frac{|S_k|}{a_k}\Big)\Big)^p(k/l)^{p\beta}\\
\leq&c_{p}C^{p}\|f\|^p(k/l)^{p\beta},
 \endaligned$$
 where
 $C$ is the same as in condition $(C_4)$ and
 $c_p$ is a positive constant
 such that for all $k$, $\dsp\ee\Big(\frac{|S_k|}{a_k}+\ee\Big(\frac{|S_k|}{a_k}\Big)\Big)^p\leq c_p$.
 Thus by the H\"older inequality, we get
 $$\aligned
 \ee\Big|\sum_{l=m}^nd_l(\xi_l-\xi_{k,l})\Big|^p
 \leq& \sum_{l_1=m}^n\cdots\sum_{l_p=m}^nd_{l_1}\cdots
 d_{l_p}(\ee|Q(l_1)|^p\cdots\ee|Q(l_p)|^p)^{1/p}\\
 \leq & c_pC^{p}\|f\|^pk^{p\beta}\sum_{l_1=m}^n\cdots\sum_{l_p=m}^nd_{l_1}\cdots
 d_{l_p}l_1^{-\beta}\cdots l_p^{-\beta}\\
 = &c_pC^{p}\|f\|^pk^{p\beta}\Big(\sum_{l=m}^nd_ll^{-\beta}\Big)^p\\
 \leq
 &c_pC^{p}\|f\|^pm^{p\beta}\Big(\sum_{l=m}^nd_l^2l\Big)^{p/2}\Big(\sum_{l=m}^nl^{-2\beta-1}\Big)^{p/2}.
 \endaligned$$
 For $m\geq 2$, it is easy to see that
\begin{align*}
m^{p\beta}\Big(\sum_{l=m}^nl^{-2\beta-1}\Big)^{p/2}
 \leq& m^{p\beta}\Big(\int_{m-1}^\infty l^{-2\beta-1}dl\Big)^{p/2}\\
 \leq&
 \Big(\frac{m}{m-1}\Big)^{p\beta}\Big(\frac{1}{2\beta}\Big)^{p/2}\leq\Big(\frac{2^{\kappa}}{\kappa}\Big)^{p/2},
\end{align*}
where $\kappa:=2\beta$. Similarly, we get for $m=1$
 $$\Big(\sum_{l=1}^nl^{-2\beta-1}\Big)^{p/2}\leq
 \Big(1+\frac{1}{\kappa}\Big)^{p/2}.$$

Hence, we have
\begin{align*}
\ee\Big|\sum_{l=m}^nd_l(\xi_l-\xi_{k,l})\Big|^p\leq
c_pC^{p}\|f\|^p\tau(\kappa)^{p/2}\Big(\sum_{l=m}^nd_l^2l\Big)^{p/2},
\end{align*}
where
$\tau(\kappa):=\Big[\Big(\frac{2^{\kappa}}{\kappa}\Big)\vee\Big(1+\frac{1}{\kappa}\Big)\Big]$.
This completes the proof of our result.
\end{proof}

\begin{lem}\label{lem4}
Assume that conditions $(C_2)-(C_4)$ of Theorem \ref{thm1} are
satisfied. Further let $b_n=0$ in condition $(C_4)$ and $f$ be a
Lipschitz function. Then for every $p\leq\mu$ and $p\in\nn$ we have
 \beq\label{230}
 \ee\Big|\sum_{k=1}^Nd_k\Big(f\Big(\frac{S_k}{a_k}\Big)-\ee
 f\Big(\frac{S_k}{a_k}\Big)\Big)\Big|^p
 \leq C_p\Big(\sum_{1\leq k\leq l\leq
 N}d_kd_l\Big(\frac{k}{l}\Big)^{\beta}\Big)^{p/2},
 \deq
 where $C_p>0$ is a constant and $\beta$ is the same as in $(C_4)$.
\end{lem}
\begin{proof} At first, we set $\dsp C_p=(4\gamma)^{p^2}$ and
$$V_{m,n}:=\sum_{l=m}^nd_ll^{-\beta}\Big(\sum_{k=1}^ld_kk^\beta\Big),\ \ (1\leq m\leq n).$$
For obtaining our result, it is enough to show that the following
claim,

"{\it if the number $\gamma$ is chosen large enough, then
 \beq\label{231}
\ee\Big|\sum_{k=m}^nd_k\xi_k\Big|^p\leq C_p(V_{m,n})^{p/2},\; \ \
{\rm for \ \ all}\ \ 1\leq m\leq n,
 \deq
where $\xi_k$ is defined as in (\ref{212}).}"

We will use induction on $p$ to show (\ref{231}). By Lemma
\ref{lem2}, we have
$$\ee\Big|\sum_{k=m}^nd_k\xi_k\Big|^2\leq 2\sum_{m\leq k\leq l\leq n}d_kd_l|\ee\xi_k\xi_l|
\leq 2c\sum_{m\leq k\leq l\leq n}d_kd_l(k/l)^\beta\leq 2cV_{m,n}.$$
Hence if we choose $\gamma$ so large that $(4\gamma)^4\geq 2c$, then
(\ref{231}) holds for $p=2$.

Assume now that (\ref{231}) is true for $p-1\geq 2$. From
$kd_k=O(1)$ it follows that there is a positive constant $A$ such
that $ \sum_{k=1}^ld_kk^\beta\geq Al^\beta$. Then we get for
$V_{m,n}\leq \gamma$ as the proof of Lemma \ref{lem3}, there exists
a constant $A_p$ such that
$$\aligned
\ee\Big|\sum_{k=m}^nd_k\xi_k\Big|^p
\leq&\sum_{k_1=m}^n\cdots\sum_{k_p=m}^nd_{k_1}\cdots
 d_{k_p}(\ee|\xi_{k_1}|^p\cdots\ee|\xi_{k_p}|^p)^{1/p}\\
 \leq& A_p\|f\|^p\sum_{k_1=m}^n\cdots\sum_{k_p=m}^nd_{k_1}\cdots
 d_{k_p}\\
 = &A_p\|f\|^p\Big(\sum_{k=m}^nd_k\Big)^p\\
 \leq
 &A_p\|f\|^pA^{-p}\Big(\sum_{k=m}^nd_kk^{-\beta}\Big(\sum_{l=1}^kd_ll^\beta\Big)\Big)^p.
 \endaligned$$
Now choose $\gamma$ so large that the
$C_p\leq({A_p}^{1/p}\|f\|/A)^p\gamma^{p/2}$. In the case of
$V_{m,n}\leq \gamma$, we have shown (\ref{231}) is valid.

 We now want to show that if for any given $X\geq \gamma$ and the inequality (\ref{231})
 holds for $V_{m,n}\leq X$, then it will also hold for $V_{m,n}\leq
 3X/2$ and this will show that (\ref{231}) holds for any value
 of $V_{m,n}$, i.e. complete the induction step.

Assume $V_{m,n}\leq 3X/2$ and set
$$S_1+S_2:=\sum_{k=m}^wd_k\xi_k+\sum_{k=w+1}^nd_k\xi_k,\ \ T_2:=\sum_{k=w+1}^nd_k\xi_{w,k},\ \ (m\leq w\leq n).$$

From the discussion of Lemma 4 in H\"ormann, S. \cite{HS} (2006),
and Lemma \ref{lem1}, for a fixed $m$ and $n$ we choose $w$ in such
a way that
$$V_{m,w}\leq X,\ \ V_{w+1,n}\leq X\ \ {\rm and}\ \ \frac{V_{w+1,n}}{V_{m,w}}=\lambda\in[1/2,1].$$

From the mean value theorem we get
 \beq\label{232}
 |S_2^j-T_2^j|\leq j|S_2-T_2|(|S_2|^{j-1}+|T_2|^{j-1})\ \ (j\geq 1).
 \deq

Since condition $(C_2)$ and Remarks \ref{rm1}, there exists a
constant $B>0$ such that for all $l\geq 1$,
$$
B\sum_{k=1}^ld_kk^\beta\geq l^{1+\beta}d_l.
$$
This also shows that
$$
 \sum_{l=m}^nld_l^2\leq BV_{m,n}, \;\;{\rm for
\;all} \;1\leq m\leq n.
$$
By Lemma \ref{lem3}, we get for all $j\geq 1$,
$$\ee|S_2-T_2|^j\leq
F_j(V_{w+1,n})^{j/2},
$$
where $F_j=B^{j/2}E_j$ and $E_j$ is the constant in Lemma
\ref{lem3}.

From the induction hypothesis in the case of $1\leq j\leq p-1$ and
from the validity of (\ref{231}) for $V_{m,n}\leq X$ in the case of
$j=p$, we have
\begin{align}
\ee|S_1|^j\leq C_j(V_{m,w})^{j/2},\;\;(1\leq j\leq p)
\end{align}
 and
 \begin{align}
 \ee|S_2|^j\leq
 C_j(V_{w+1,n})^{j/2}\leq C_j\lambda^{j/2}(V_{m,w})^{j/2},\;\;(1\leq j\leq p).
\end{align}
The remains of the proof are the same as in Lemma 4 in H\"ormann, S.
\cite{HS} (2006), but for completeness, we still give the proof.  By
$C_r$ inequality, we have
 \beq\label{234}
 \ee|T_2|^j\leq2^jC_j\lambda^{j/2}(V_{m,w})^{j/2},\;\; (1\leq j\leq
 p).
 \deq
Furthermore, from H\"older inequality the following result is easy,
 \begin{eqnarray}\label{235}
 \ee|S_1|^j|S_2-T_2||S_2|^{p-j-1}&\leq&(\ee|S_1|^p)^{j/p}(\ee|S_2-T_2|^p)^{1/p}(\ee|S_2|^p)^{(p-j-1)/p}\nonumber\\
 &\leq&C_p^{(p-1)/p}F^{1/p}_p\lambda^{(p-j)/2}(V_{m,w})^{p/2}.
 \end{eqnarray}
The last inequality remains valid, with an extra factor $2^{p-j-1}$
on the right hand side, if $|S_2|^{p-j-1}$ on the left hand side is
replaced by $|T_2|^{p-j-1}$. Since $S_1$ and $T_2$ are independent,
we get
$$\ee|S_1+S_2|^p
\leq
\ee|S_1|^p+\ee|S_2|^p+\sum_{j=1}^{p-1}G_p^j(\ee|S_1|^j|S_2^{p-j}-T_2^{p-j}|+\ee|S_1|^j\ee|T_2|^{p-j}),$$
where $G_p^j$ denote the combination, i.e.,
$G_p^j=p![j!(p-j)!]^{-1}$. Substituting $(\ref{232})-(\ref{235})$
(using also the analogue of $(\ref{235})$ with $|T_2|^{p-j-1}$) in
the above inequality and get
$$\aligned
\ee|S_1+S_2|^p \leq&
C_p(V_{m,w})^{p/2}\Big(1+\lambda^{p/2}+C_p^{-1/p}F^{1/p}_p
\sum^{p-1}_{j=1}2^{p-j}G_p^j(p-j)\lambda^{(p-j)/2}\\
&+C^{-1}_p\sum^{p-1}_{j=1}2^{p-j}G_p^j\lambda^{(p-j)/2}C_jC_{p-j}\Big).
\endaligned$$
Note that
$$
C_p^{-1/p}F^{1/p}_p\leq const\cdot
\tau(\kappa)^{1/2}c_p^{1/p}(4\gamma)^{-p}
$$
and
$$
C_jC_{p-j}/C_p\leq(4\gamma)^{-p}, \;(1\leq j\leq p-1),
$$
thus, by $\lambda\leq 1$, we have
$$
C_p^{-1/p}F^{1/p}_p
\sum^{p-1}_{j=1}2^{p-j}G_p^j(p-j)\lambda^{(p-j)/2}\leq const\cdot
\tau(\kappa)^{1/2}p\,c_p^{1/p}\gamma^{-p}
$$
and
$$C^{-1}_p\sum^{p-1}_{j=1}2^{p-j}G_p^j\lambda^{(p-j)/2}C_jC_{p-j}\leq const\cdot
\gamma^{-p}.$$ Since $\lambda\geq 1/2$ we have shown that for a
large $\gamma$ the relation $\ee|S_1+S_2|^p\leq
C_p(1+\lambda)^{p/2}(V_{m,w})^{p/2}=C_p(V_{m,n})^{p/2}$ is true,
i.e., for $V_{m,n}\leq 3X/2$, (\ref{231}) is valid.
\end{proof}
\begin{lem}\label{lem5}(See Lemma 5 in \cite{HS})
Assume the condition $(C_3)$ of Theorem \ref{thm1} is satisfied.
Then for any $\alpha>0$ and any $\eta<\rho$, we have
$$\sum_{1\leq k\leq l\leq N}d_kd_l\Big(\frac{k}{l}\Big)^\alpha=O\Big(\frac{D_N^2}{(\log D_N)^\eta}\Big).$$
\end{lem}

\begin{proof} [{\bf Proof of Theorem \ref{thm1}}]
Without loss of generality, from Lemma \ref{lem4} and Lemma
\ref{lem5}, we have, for any $\varepsilon>0$, $p\leq\mu$ and
$p\in\nn$,
$$\pp\Big(\Big|\sum_{k=1}^Nd_k\xi_k\Big|>\varepsilon D_N\Big)\leq c(p,\varepsilon)(\log D_N)^{-p\eta/2},
\;\ \ for \ \ N\geq N_0.
$$
Since $\mu\geq(2\vee 4/(\rho-r))$ for some $0<r<\rho$, we can choose
suitable $\eta<\rho$ and $p$ such that $p\eta> 4$. By $(C_3)$, we
have $D_{N+1}/D_N\to 1$, thus we can choose $(N_j)$ such that
$D_{N_j}\sim\exp\{\sqrt{j}\}$. Applying Borel-Cantelli lemma, we get
$$
\lim_{j\to\infty}D_{N_j}^{-1}\sum_{k=1}^{N_j}d_k\xi_k=0\;\ \ a.s..
$$
For $N_j\leq N\leq N_{j+1}$, we have
$$
D_N^{-1}|\sum_{k=1}^{N}d_k\xi_k|\leq
D_{N_j}^{-1}|\sum_{k=1}^{N_j}d_k\xi_k|+2(D_{N_{j+1}}/D_{N_j}-1)\;\ \
a.s..
$$
The convergence of the subsequence implies that whole sequence
converges a.s., since $D_{N_{j+1}}/D_{N_j}\to 1$. This complete the
proof of the theorem.
\end{proof}

{}

\end{document}